\documentclass[twocolumn]{revtex4}

\usepackage{epsfig}
\usepackage{graphicx}
\usepackage{amsmath}
\usepackage{mathrsfs}
\usepackage{dcolumn}
\usepackage{amsfonts}
\usepackage{ulem}
\usepackage{booktabs}
\usepackage{makecell}

\newcommand{\ket}[1]{|#1\rangle}
\newcommand{\bra}[1]{\langle#1|}

\begin{document}
 \title{Quantum Lyapunov control with machine learning}
 \author{S. C. Hou }
 \author{X. X. Yi }
 \email{ yixx@nenu.edu.cn}

\affiliation{Center for Quantum Sciences and School of
Physics, Northeast Normal University, Changchun 130024, China}

\begin{abstract}
Quantum state engineering is a central task in Lyapunov-based quantum control.
Given different initial states, better performance may be achieved if the control
parameters, such as the Lyapunov function, are individually optimized for each
initial state, however, at the expense of computing resources. To tackle this issue, we propose an initial-state-adaptive Lyapunov control strategy with machine learning, specifically, artificial neural networks trained through supervised learning.  Two designs are presented and illustrated where the feedforward neural network and the general regression neural network are used to select control schemes and design Lyapunov functions, respectively. Since the sample generation and the training of neural networks are carried out in advance, the initial-state-adaptive Lyapunov control can be implemented without much increase of computational resources.
\end{abstract}
\date{\today}
\maketitle

\section{Introduction}

Quantum control \cite{Alessandro2007,Wiseman2009} plays a fundamental role in modern quantum technologies such as quantum computation, quantum communication and quantum metrology. A central goal in quantum control is designing time-varying external control fields to effectively engineer quantum states and operators. More than a decade ago, the Lyapunov-based method was developed for the control of quantum systems  \cite{Vettori2002,Grivopoulos2003}. In quantum Lyapunov control, the control fields are obtained
by simulating the dynamics (in feedback form) only once and then applied in an open-loop scenario.
The method has the merits of simplicity in generating control fields and flexibility in designing the control field shapes. In recent years, numerous efforts have been devoted to investigate or improve the convergence of Lyapunov control for different quantum control models \cite{Mirrahimi2005, Kuang2008, Wang2010, Hou2012, Wang2014, Zhao2012, Silveira2016, Kuang2017}. Meanwhile,  Lyapunov control method is successfully employed for diverse quantum information processing tasks \cite{ Wang2009, Yi2009, Sayrin2011, Dong2012, Hou2014, Shi2015PRA, Shi2015SR, Shi2016, Li2016, Silveira2016, Ran2017,Li2018}. For example, it is recently used to realize topological modes \cite{Shi2015SR},  quantum synchronization \cite{Li2016} and speed up adiabatic passage \cite{Ran2017}.

In previous research, the designs of quantum Lyapunov control are usually initial-state-independent.  When dealing with different initial states, better performance (e.g. control time or fidelity) may be achieved if the control parameters (such as those in the Lyapunov function) are individually optimized for each initial state. However, it is usually hard to find an explicit (analytic) relationship between the optimal control parameters and the initial states since the control fields are generated numerically. On the other hand, numerical optimizing the Lyapunov control parameters  typically requires simulating the dynamics more than once, making quantum Lyapunov control complicated. Thus an initial-state-adaptive quantum Lyapunov control without a significant increase of computing resources is desirable.

Machine learning \cite{Haykin2009, Alpaydin2010} is a powerful tool to improve a performance criterion from experience or data, which has been extensively applied in internet technology, artificial intelligence, finance, medical diagnosis and so on. Recently, machine learning technology has been successfully employed to advance quantum physics problems \cite{Magesan2015, Mills2017, Melnikov2018, Torlai2018, Carleo2017, Deng2018, Gao2018, Zahedinejad2016, Mavadia2017, August2017, Yang2018}, such as quantum many-body problems \cite{Carleo2017, Deng2018}, quantum state identification \cite{Deng2018,Gao2018} and quantum control \cite{Zahedinejad2016, Mavadia2017, August2017, Yang2018}. Motivated by its ability and versatility, we intend to use machine leaning techniques, specifically, (artificial) neural networks \cite{Haykin2009} to design an initial-state-adaptive quantum Lyapunov control. The basic idea is as follows. First, numerically generate a certain number of samples that encode different initial states and their corresponding optimal parameters. Then, train a neural network with these samples through supervised learning until its performance is satisfactory. At last, apply the trained neural network to predict control parameters for new initial states. Two designs are proposed to select control schemes and design Lyapunov functions. The initial-state-adaptive designs would be helpful when the number of initial states is large or real-time control is needed.

The remainder of the paper is organized as follows.  Sec.{\rm II} reviews the Lyapunov control method for the eigenstate preparation problem. In Sec.{\rm III}, we introduce the feedforward neural network (multilayer perceptron) and the general regress neural network (GRNN) which are used as the tools for classification and regression in this paper, respectively. The two initial-state-adaptive designs pare proposed in Sec.{\rm IV} and then illustrated with a three-level quantum system in Sec.{\rm V}. Finally, the results are summarized and discussed in Sec.{\rm VI}.

\section{quantum Lyapunov control }
Quantum Lyapunov control is a useful technique for quantum control tasks, typically eigenstate control \cite{Wang2009, Dong2012, Sayrin2011, Shi2015PRA, Shi2016, Ran2017, Li2018}.  It consists of two steps. In the first step, time-dependent control fields are numerically calculated by simply one simulation of the system dynamics (in feedback form). In the second step, the generated control fields are used in applications (experiments) in an open-loop way.

We introduce the mathematical formula of quantum Lyapunov control with a n-dimensional closed quantum system described by the Schr\"{o}dinger equation ($\hbar=1$ is assumed)
\begin{eqnarray}
\frac{d}{dt}\ket{\Psi}=-i[H_0+\sum_{k=1}^{m}f_k(t)H_k]\ket{\Psi}.
\label{eqn:LyaCtrlEqn}
\end{eqnarray}
Here $H_0$ is the system (drift) Hamiltonian, $H_k$ is the $k$th control Hamiltonians and $f_k(t)$ is its corresponding control field which is a time-dependent real function. The aim is to find proper $f_k(t)$ such that the initial state $\ket{\Psi}$ evolves to a desired state  $\ket{\Psi_d}$  at some point of time.

In quantum Lyapunov control, a real function $V$ called Lyapunov function (conventionally $V \geq 0$) is assigned and $f(t)$ are designed to guarantee $\dot{V}\leq0$. Through this, the quantum system is driven to states satisfying $\dot{V}=0$ as $t\rightarrow\infty$, meanwhile, the desired state is asymptotically reached. The convergence behavior could be analyzed by the La Salle's invariance principle \cite{Alessandro2007}. The choice of Lyapunov function $V$ is not unique. For example, $V$ could be chosen as the distance between the quantum state and the desired state, the expectation value of a Hermitian operator and so on \cite{Kuang2008}. Here we consider the second form of Lyapunov function, i.e.,
\begin{eqnarray}
V=\bra{\psi} P \ket{\psi},
\label{eqn:LyaFun}
\end{eqnarray}
where $P$ is a Hermitian and positive semi-definite operator such that $V \geq 0$. This form is representative and  covers some other forms of Lyapunov function such as that based on the Hilbert-Schmidt distance \cite{Alessandro2007}. More importantly, there is freedom in designing $P$ enabling us to optimize it for different initial states for the purpose of this paper.

The control fields could be designed based on the time derivative of $V$,
\begin{eqnarray}
\dot{V}&=&\bra{\Psi}i[H_0+\sum_{k=1}^{m}f_k(t)H_k,P]\ket{\Psi}\\
&=&\sum_{k=1}^m f_k(t)\bra{\Psi}i[H_k,P]\ket{\Psi}
\label{eqn:dV}
\end{eqnarray}
where $[H_0,P]=0$ is assumed to cancel the drift term. This condition could be realized by constructing the hermitian operator $P$ as
\begin{eqnarray}
P=\sum_{l=1}^{n}p_l\ket{E_l}\bra{E_l}.
\label{eqn:P}
\end{eqnarray}
Here $\ket{E_l}$ is the $l$th eigenstate of $H_0$ and $p_l$ are non-negative real numbers.
In this work, $p_l$ will be optimized for different initial states and predicted by trained artificial neural networks.

The control fields $f_k(t)$ is conventionally designed as
\begin{eqnarray}
f_k(t)= -K \bra{\Psi}i[H_k,P]\ket{\Psi}
\label{eqn:ControlField}
\end{eqnarray}
where $K$ is a real constant associated with the control strength. Other approaches to design the control fields are also investigated to improve the performance of quantum Lyapunov control \cite{Hou2012, Zhao2012, Kuang2017}.
From Eq(\ref{eqn:ControlField}), there is
\begin{eqnarray}
\dot{V}=-K^{-1}\sum_{k=1}^m f_k^2(t)\leq0,
\label{eqn:dVless0}
\end{eqnarray}
i.e., the Lyapunov function keeps non-increasing with the controlled dynamics.
With ideal control parameters (e.g. Lyapunov function, control Hamiltonian, design of control fields), the control law determined by Eq.(\ref{eqn:LyaFun},\ref{eqn:P},\ref{eqn:ControlField}) will drive any initial state $\ket{\Psi(0)}$ (except that satisfies $\dot{V}(0)=0$) asymptotically to the eigenstate of $H_0$ with the minimum eigenvalue as $t\rightarrow\infty$. Meanwhile, $V$ will decrease to its minimum.

Obviously, the performance (e.g., fidelity, control time) of quantum Lyapunov control depends on the control parameters such as Lyapunov function $V$ and control Hamiltonian $H_k$. Choosing appropriate parameters is therefore of great importance for Lyapunov control problems.

\section{artificial neural networks}
In this section, we briefly introduce two neural network models used in this paper, feedforward neural network and general regress neural network. Mathematically, these neural networks could be understood as a function that maps an input real vector $X$ to an output real vector $Y$.

\subsection{Feedforward Neural Network}
\begin{figure}
\includegraphics*[width=8.5cm]{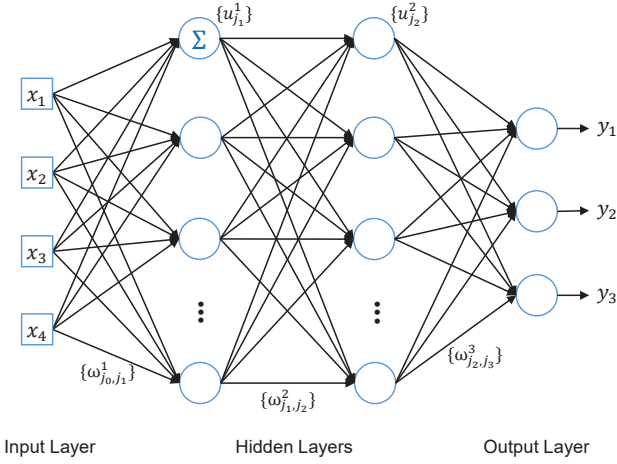}
\caption{A schematic diagram of a feedforward neural network 4 input nodes, 3 output neurons and 2 hidden layers. Here $j_l=1,2,...,n_l$ where $n_l$ is the node (neuron) number of the $l$th layer.}
\label{FIG:MPNN}
\end{figure}

Feedforward neural network is the most well known neural network. A schematic diagram of a feedforward neural network is shown in Fig.\ref{FIG:MPNN}. A feedforward neural network consists of a layer of input nodes (by squares in Fig.\ref{FIG:MPNN}), an output layer of neurons (processing units, by circles in Fig.\ref{FIG:MPNN}), and possibly a set of hidden layers of neurons. In feedforward neural networks, signal flows from the input layer to the output layer without feedback loops. A feedforward neural network with one or more hidden layers is called a multilayer perceptron \cite{Alpaydin2010,Haykin2009}. With enough neurons, a multilayer perceptron is able to approximate any continuous nonlinear function and solve many complicated tasks. In a feedforward neural network, the output $y$  of a single neuron is expressed by
\begin{eqnarray}
y=s(\sum_{i=1}^{m}x_{i}\omega_i +b),
\label{eqn:neuron}
\end{eqnarray}
where the $x_i$ is the output of the $i$th neuron (node) of the last layer, $\omega_i$ is the weight of $x_i$ corresponding to the arrows in Fig.\ref{FIG:MPNN}, and $b$ is a bias (threshold) which is omitted in Fig.\ref{FIG:MPNN}. $s(...)$ is called an activation function which is usually a nonlinear sigmoid function limiting the strength of the output signal. The logistic function
\begin{eqnarray}
s(x)=\frac{1}{1+e^{-x}}
\label{eqn:sigmoid}
\end{eqnarray}
is used as the activation function in this paper that transfers any input signal to the range $0$ to $1$ .

In a feedforward neural network with $m$ layers of neurons ($n_l$ neurons in the $l$th layer) and  an input layer with $n_0$ nodes, the input $X=[x_1\ x_2\ ... \ x_{n_0}]^T$ is transformed to the output $Y=[y_1\ y_2\ ... \ y_{n_m}]^T$  by
\begin{eqnarray}
u^1_{j_1}&=&s(\sum_{j_0=1}^{n_0}x_{j_0}  w_{j_0,j_1}^{1}+b^1_{j_1}), \\ u^2_{j_2}&=&s(\sum_{j_1=1}^{n_1}u^1_{j_1} w_{j_1,j_2}^{2}+b^2_{j_2}), \\
&\vdots& \nonumber\\
y_{j_m}&=&s(\sum_{j_{m\!-\!1}=1}^{n_{m\!-\!1}}   u^m_{j_{m\!-\!1}} w_{j_{m\!-\!1},j_m}^{m}+b^m_{j_m})
\label{eqn:NN}
\end{eqnarray}
where $j_l=1,2,...,n_l$ with $l=0,1,...,m$. Here $u^\alpha_{j_l}$ is the output of the $j_l$th neuron of the $\alpha$th neuron layer. In the above equations, the superscript (1,2,...,m) denotes the index of the neuron layers, and the subscript $j_l$ represents the $j_l$th neuron or node of the $l$th layer, as shown in Fig.\ref{FIG:MPNN}. Thus the network is determined by the layer numbers m, node numbers in each layers, weights, biases as well as the activation function.

For a specific problem, the design of the feedforward neuron network structure is generally empirical. The training of a feedforward neuron network is implemented by adjusting its weights and biases. In supervised learning, the weights and biases could be effectively studied by the back-propagation (BP) algorithm \cite{Haykin2009, Alpaydin2010} with the training samples  including a number of input vectors and their target output vectors.

\subsection{General Regression Neural Network}
\begin{figure}
\includegraphics*[width=8.5cm]{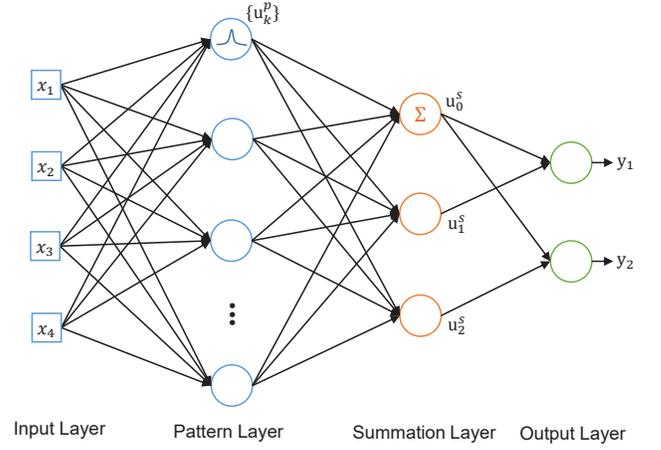}
\caption{A schematic diagram of a general regression neural network with 4 input nodes and  2 output neurons.}
\label{FIG:GRNN}
\end{figure}
General regression neural network (GRNN) is a type of radial basis function (RBF) network proposed by  D. F. Specht in 1991 \cite{Specht1991}. It is a powerful tool to estimate continuous variables \cite{Leung2000, Li2011, Liu2014, Panda2015} even when the training data is few. A general regression neural network consists of 4 layers: an input layer, a pattern layer, a summation layer and an output layer as shown in Fig.\ref{FIG:GRNN}. In contrast to the feedforward neural network, the neuron number in each layer is fixed in a GRNN and determined from its training samples. For a training set with $N$ samples $\{X^k,Y^k\}, k=1,2,...,N$ where $X^k=[x^k_1\  x^k_2 \ ...\  x^k_{n_I}]$ is the $k$th input vector and $Y^k=[y^k_1\  y^k_2 \ ...\  y^k_{n_O}]$ is its target output, the pattern layer neuron number is $N$, the number of the input (output) layer node is $n_I$ ($n_O$), and the number of the summation layer neuron is $n_O+1$.

For an input $X$, the output of the pattern layer ($u^p_k$, $k=1,2,...,N$), the summation layer ($u^s_{j_s}$, $j_s=0,1,2,...,n_O$) and the output layer ($y_j$, $j=1,2,...,n_O$) is given by
\begin{eqnarray}
u^p_{k}=\exp(-\frac{(X-X^k)^T(X-X^k)}{2\sigma^2}),\label{eqn:GRNN1}\\
u^{s}_0=\sum_{k=1}^{N} u^p_k, \quad u^s_{j_s (j_s\neq0)}=\sum_{k=1}^{N} y^k_{j_s} u^p_k, \label{eqn:GRNN2}\\
y_j=\frac{ u^s_j }{u^s_0},\quad\quad\quad\quad\quad\quad\quad\label{eqn:GRNN3}
\end{eqnarray}
respectively. The estimation $y_j$ for an input $X$ can be understood as an average of all $y^k_j$ weighted exponentially according to the Euclidean distance between $X$ and $X^k$. Here $\sigma$ is called a smoothing parameter ($\sigma>0$). When $\sigma$ is small, the estimation $y_i$ for $X$ is closely related to $y^k_i$ whose inputs $X^k$ are close to $X$.  In contrast, when $\sigma$ is large, $y_i$ approaches the mean of all  $y^k_j$. The GRNN is established as soon as the training samples are stored while the smoothing parameter $\sigma$ is the only adjustable parameter needed to train.

\section{Initial-state-adaptive designs }
When dealing with different initial states, better performances may be achieved if the parameters in quantum Lyapunov control, such as those in the Lyapunov function $V$ or the control Hamiltonian $H_k$, are optimally chosen for each initial state, i.e., initial-state-adaptive. However, finding optimized parameters typically costs more computing resources (e.g. simulation time) since Lyapunov control fields are calculated numerically. In this way, quantum Lyapunov control would lose its simplicity that the control fields are obtained with only one simulation of system dynamics. To tackle this issue, we propose to design an initial-state-adaptive Lyapunov control with neural networks. As the processes of generating samples and training neural networks are implemented in advance, the computing resources of initial-state-adaptive Lyapunov control with neural networks would not significantly increase in applications. The basic strategy comprises the following steps.\\
(1) Generate a certain number of initial states whose parameters are randomly distributed in an interested ranges.\\
(2) Numerically find the optimal control parameters for these initial states.  Data from (1) and (2) constitute the training set (testing set if needed).\\
(3) Build a neural network whose input is associated with the initial state parameters and output associated with the optimal control parameters.\\
(4) Train the neural network by supervised learning with the training set until the neural network performance is satisfactory.\\
(5) Apply the trained neural network to predict optimal control parameters for new initial states .\\

For the eigenstate control problem described in Sec.II, we propose two designs where the mentioned neural networks are used to select control schemes or predict Lyapunov functions for different initial states. Two important functions of neural networks, classification and regression, are employed.

\subsection{Classification: selecting control schemes}
Consider there are several Lyapunov control schemes where the control Hamiltonians, Lyapunov functions or other conditions are different. One of these schemes will be finally adopted in experimental or theoretical application. Our aim is to use neural networks to predict the optimal scheme for each individual initial state.

Specifically, assume there are $M$ candidate control Hamiltonians $H_c (c=1,2...M)$ in which one of them will be selected. Thus the dynamics is described by
 \begin{eqnarray}
\frac{d}{dt}\ket{\Psi}=-i[H_0+f(t)H_c]\ket{\Psi}.
\label{eqn:OneCtrl}
\end{eqnarray}
Other conditions such as the Lyapunov function $V$, the strength $K$, and $H_0$ are fixed.
The task is preparing an eigenstate of $H_0$, say, $\ket{E_g}$, as discussed in Sec.II. Given an initial state, we will use a feedforward neural network to predict the control Hamiltonian that leads to the highest fidelity defined by $F=|\bra{\Psi(T)}E_g\rangle|^2$ at a certain control time $T$. The problem is solved by classifying the initial states according to their favorable control Hamiltonian with the feedforward neural network.

For a $n$-dimensional system,  the initial state $\ket{\Psi_0}$ could be parameterized in the eigenbasis of $H_0$ as
 \begin{eqnarray}
 \begin{split}
\ket{\Psi_0}_{n=2} =&\sin\theta_1e^{i\phi_1}\ket{E_1} +\cos\theta_1\ket{E_2} \\
\ket{\Psi_0}_{n=3} =&\sin\theta_2(\sin\theta_1e^{i\phi_1}\ket{E_1}+\cos\theta_1e^{i\phi_2}\ket{E_2})\\
 &+\cos\theta_2\ket{E_3}\\
\vdots\ &
 \end{split}
\label{eqn:initial}
\end{eqnarray}
where $0\leq\theta_i\leq\frac{\pi}{2}$ and  $0\leq\phi_i\leq 2\pi$. It is observed that $2(n-1)$ parameters, $\theta_{1,2,...,n-1}$ and $\phi_{1,2,...,n-1}$,  are required to
determine an initial state up to a non-physical global phase.

We define the training set with $N_{train}$ samples as
\begin{eqnarray}
S=\{(X^1,Y^1),(X^2,Y^2),...,(X^{N_{train}},Y^{N_{train}})\}.
\label{eqn:TraingSet}
\end{eqnarray}
In the $k$th sample, $X^k$ is the input vector with $2(n-1)$ elements defined as
 \begin{eqnarray}
X^k
=[\theta^k_{1}\ \theta^k_{2}\ ...\  \theta^k_{(n-1)}\  \phi^k_{1}\ \phi^k_{2}\ ... \ \phi^k_{(n-1)} ]^T.
\label{eqn:input}
\end{eqnarray}
In this paper, we assume all the possible initial states are interested, thus $\theta^k_i$ and $\phi^k_i$ could be chosen as random numbers uniformly distributed in $[0,\frac{\pi}{2}]$ and $[0,2\pi]$, respectively. For an initial state, the choice of its favorable control Hamiltonian is determined by simulating the dynamics with $M$ candidate control Hamiltonians and comparing the fidelities. The target output vector $Y^k$ indicating the choices is a unit vector $\vec{e}$ with $M+1$ elements.  For example, the control Hamiltonian could be mapped to the output vector $Y$ by
\begin{eqnarray}
\begin{split}
H_1 \longrightarrow \vec{e}_1&=&[1\  0\ ... \ 0]^T, \\
 H_2 \longrightarrow \vec{e}_2&=&[0\  1\  ... \ 0]^T, \\
&\vdots&   \\
 H_M \longrightarrow \vec{e}_M&=&[0\  ... \ 1\  0]^T, \\
  others  \longrightarrow  \vec{e}_{M+1} &=&[0\  ... \ 0\  1]^T,
\end{split}
\label{eqn:outputmap}
\end{eqnarray}
where \textit{others} refer to cases without an optimal choice, e.g., ineffective controls or the existence of  the same fidelities. On the other hand, a testing set $S_T$ with $N_{test}$ samples could be generated in a similar way as the training set eq.(\ref{eqn:TraingSet}). The testing set (not participate in the supervised learning) is used for checking the performance of the neural network in order to avoid over training of the neural network.

Next, a feedforward neural network with $2(n-1)$ input nodes, $M+1$ output neurons, plus some hidden layers is set up with the activation function Eq.(\ref{eqn:sigmoid}).  For an input vector $X$, the output of the neural network is a linear combination of all the basis vectors, i.e.,
\begin{eqnarray}
\begin{split}
Y'=\sum_{j=1}^{M+1} q_j \vec{e}_k, \quad 0<q_j<1.
\end{split}
\label{eqn:map}
\end{eqnarray}
The classification is implemented by selecting the choice with the largest coefficient $q_j$. Here $q_j$ might be understood as an unnormalized probability that the choice is $j$. The performance of a neural network could be measured by the mean squared error (MSE)
\begin{eqnarray}
MSE=\frac{1}{N}\sum_{k=1}^{N}(Y'^k-Y^k)^T(Y'^k-Y^k)
\label{eqn:MSE}
\end{eqnarray}
where $Y'^k$ is the output of the neural network for $X^k$ and $Y^k$ is the $k$th target output vector. $N$ is the number of the training (or testing) samples.

With the training set Eq.(\ref{eqn:TraingSet}) determined by Eq.(\ref{eqn:input}) and Eq.(\ref{eqn:outputmap}), the weights and biases could be effectively trained by the back propagation (BP) algorithm. The iteration number of the BP training process could be determined by checking the mean squared error (or the classification success rate) for the testing set. Before the training (testing) process, the input vector $X^k=[x_1^k\ x_2^k\ ...\ x_{n_I}^k]$ is normalized to $X'^k=[x'^k_1\ x'^k_2\ ...\ x'^k_{n_I}]$ by $x'^k_j=2(x^k_j-x_j^{min})/(x_j^{max}-x_j^{min})-1$ where $x_j^{max}$ ($x_j^{min}$) is the maximum (minimum) of the $N_{train}$ input vector elements $x^k_j$ ($k=1,2,...,N_{train}$) from the training set. In this way, all the signals sent to the input nodes are scaled to the range $[-1,1]$ in order to be sensitive to the sigmoid functions of the neural network. Finally, the trained neural network will be used to select control Hamiltonian for new initial states (out of the training set). For the problem of selecting other control schemes, our method may also be applied in a similar way.

\subsection{Regression: designing Lyapunov function}

In this section, GRNN is used to design an initial-state-adaptive Lyapunov function $V$  of the form Eq.(\ref{eqn:LyaFun}) where $P=f_{GRNN}(\ket{\Psi_0})$. The system Hamiltonian $H_0$ and control Hamiltonian(s) $H_k$ are fixed.  The task is to prepare an eigenstate of $H_0$ with a high fidelity defined as $F=|\bra{\Psi(T)}E_g\rangle|^2$ at time $T$.

Notice that the strength coefficient $K$ in Eq.(\ref{eqn:ControlField}) can be absorbed into the operator $P$, i.e., $V'=\bra{\psi} KP \ket{\psi}=\bra{\psi} P' \ket{\psi}$. Therefore, we set $K=1$ and merely discuss $P$ for simplicity. Assume the goal state is the $g$th eigensate of $H_0$ denoted by $\ket{E_g}$. For a $n$-dimensional system, the operator $P$ is designed as
\begin{eqnarray}
\begin{split}
P &=& p_g\ket{E_g}\bra{E_g}+\sum_{l\neq g}p_l\ket{E_l}\bra{E_l}\quad\quad\\
  &=& \sum_{l\neq g}p_l\ket{E_l}\bra{E_l} \quad \quad(p_l>p_g=0).
\end{split}
\label{eqn:Pdesign}
\end{eqnarray}
We have set the minimum coefficient $p_g$ to $0$ without loss of generality, since if $p_g\neq0$, one can shift it to zero by adding $-p_g\sum_{l}\ket{E_l}\bra{E_l}=-p_g\text{I}$ to $P$, which does not change the control fields according to Eq.(\ref{eqn:ControlField}). Now the favorable Lyapunov function for an initial state could be obtained by optimizing  $p_{l(l\neq g)}$ numerically. The number of $p_l$ to be optimized is $n-1$. In principle, there are no limitation on the bound of $p_l$ in optimization. However, some constraints of $p_l$ are required to limit the strength of the control fields and facilitate the numerical optimizations, e.g., $0<p_l\leq p_l^{max}$.

In this method, the training set (testing set if needed) with $N_{train}$ ($N_{test}$) samples can be defined similarly as  Eq.(\ref{eqn:TraingSet}). The input vector $X^k$ is given by Eq.(\ref{eqn:input}) and the $k$th output vector $Y^k$ is defined as
 \begin{eqnarray}
Y^k=[p^k_1\ p^k_2\ ... \ p^k_{n-1}]^T.
\label{eqn:Voutput}
\end{eqnarray}
The elements of the target vector $p^k_{j=1,2,...,n-1}$ are the optimal values in one-to-one correspondence with $p_l$ in Eq.(\ref{eqn:Pdesign}). Given an input vector $X^k$  (an initial state), $Y^k$ is obtained by numerically finding $p^k_j$ that maximizing the fidelity $F_k$ with the constraints of $p^k_j$, in which the dynamics needs to be simulated many times.

With the training set, a GRNN  with $2(n-1)$ input nodes, $N_{train}$ pattern layer neurons and $n-1$ output neurons can be built straightforwardly. The smoothing parameter $\sigma$ (defined in Eq.(\ref{eqn:GRNN1})) could be determined by checking the GRNN performance for the testing samples without many trials. Different performances for the testing samples might be used such as the MSE or the averaged logarithmic infidelity defined by
 \begin{eqnarray}
\epsilon=\frac{1}{N}\sum_{k=1}^{N}\log(1-F_k).
\label{eqn:AveLogInf}
\end{eqnarray}
Since the smoothing parameter $\sigma$ (connected with the width of the Gaussian function in Eq.(\ref{eqn:GRNN1})) is the same for all the input vector elements, the input vectors need to be normalized such that $\sigma$ is sensitive to
all the input vector elements. Let each input vector element be normalized to $[-1,1]$ as in the last design, then the average space between two neighboring (normalized) input vectors could be estimated by $D=2/\sqrt[n_I]{N_{train}}$ where $n_I=2(n-1)$ is the dimension of the input vector.  We suggest to find the smoothing parameter in the range $0<\sigma< \sigma_{max} $ where $\sigma_{max}\sim D$. With an appropriate $\sigma$, the GRNN will be finally used to predict optimal $P$ of the Lyapunov function for new initial states.


\section{Illustration}
In this section, we illustrate our designs with a three-level quantum system (n=3). The time-independent system Hamiltonian is described by
\begin{eqnarray}
H_0=\sum_{n=1}^3\omega_n\ket{n}\bra{n}+g(\ket{1}\bra{2}+\ket{2}\bra{1}),
\label{eqn:H0example}
\end{eqnarray}
where $\omega_n$ is the frequency of the $n$th level  and the state $\ket{2}$ and $\ket{1}$ are coupled with a strength $g$. The task is to prepare an eigenstate of $H_0$ with high fidelity, e.g., $\ket{E_3}$ where $E_3$ is the highest eigenenergy. The fidelity is calculated by $F=|\bra{\Psi(T)}{E_3}\rangle|^2$ with $T$  the control time. The dynamics is described by Eq.(\ref{eqn:OneCtrl}) and the control field is given by Eq.(\ref{eqn:ControlField}).

\subsection{Selecting control Hamiltonians}
To illustrate the first design, consider two candidate control Hamiltonians ($M=2$),
\begin{eqnarray}
H_1&=&\ket{1}\bra{3}+\ket{3}\bra{1},\label{eqn:Hami1}\\
H_2&=&\ket{2}\bra{3}+\ket{3}\bra{2}.\label{eqn:Hami2}
\end{eqnarray}
Given an arbitrary initial state $\ket{\Psi_0}$, we will use a feedforward neural network to select the control Hamiltonian that leads to a higher fidelity. In this example, $P=\ket{E_1}\bra{E_1}+\ket{E_2}\bra{E_2}=\text{I}-\ket{E_3}\bra{E_3}$ is used in the Lyapunov function of the form Eq.(\ref{eqn:LyaFun}). Thus $V=\bra{\Psi} P \ket{\psi}=1-|\bra{\Psi}E_3\rangle|^2$  might be understood as either Eq.(\ref{eqn:Pdesign}) with $p_1=p_2=1$ (unoptimized) or a distance between the controlled state and the goal state. According to our method, the input vectors in the training (testing) set are given by $X^k=[\theta^k_{1}\ \theta^k_{2}\ \phi^k_{1}\ \phi^k_{2}]^T$, $k=1,2,...,N_{train} (N_{test})$ where $\theta^k_{1,2}$ and $\phi^k_{1,2}$ is defined as in Eq.(\ref{eqn:initial}) with $n=3$. We consider 3 choices ($H_1$, $H_2$ and \textit{others}) corresponding to target output vectors $\vec{e}_{1,2,3}$, respectively. Here the \textit{others}) refers to the same fidelities or both fidelities less than $0.99$.

\begin{figure}
\includegraphics*[width=8cm]{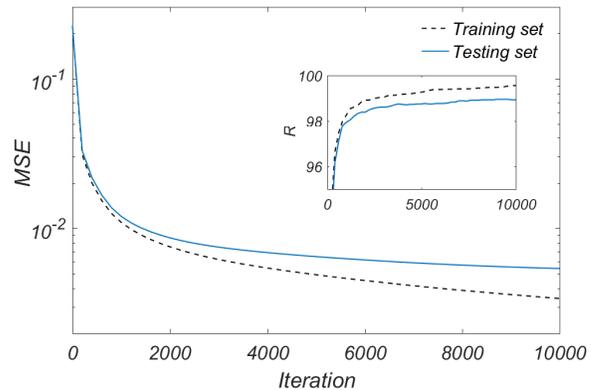}
\caption{The training process of the feedforward neural network. Mean figure: the Mean squared error (MSE) for the training set (black dashed line) and the testing set (blue solid line) versus the iteration number. Inset: the percentage of the classification success rate (R) for the training set and the testing set. }
\label{FIG:Training}
\end{figure}

\begin{figure}
\includegraphics*[width=4.2cm]{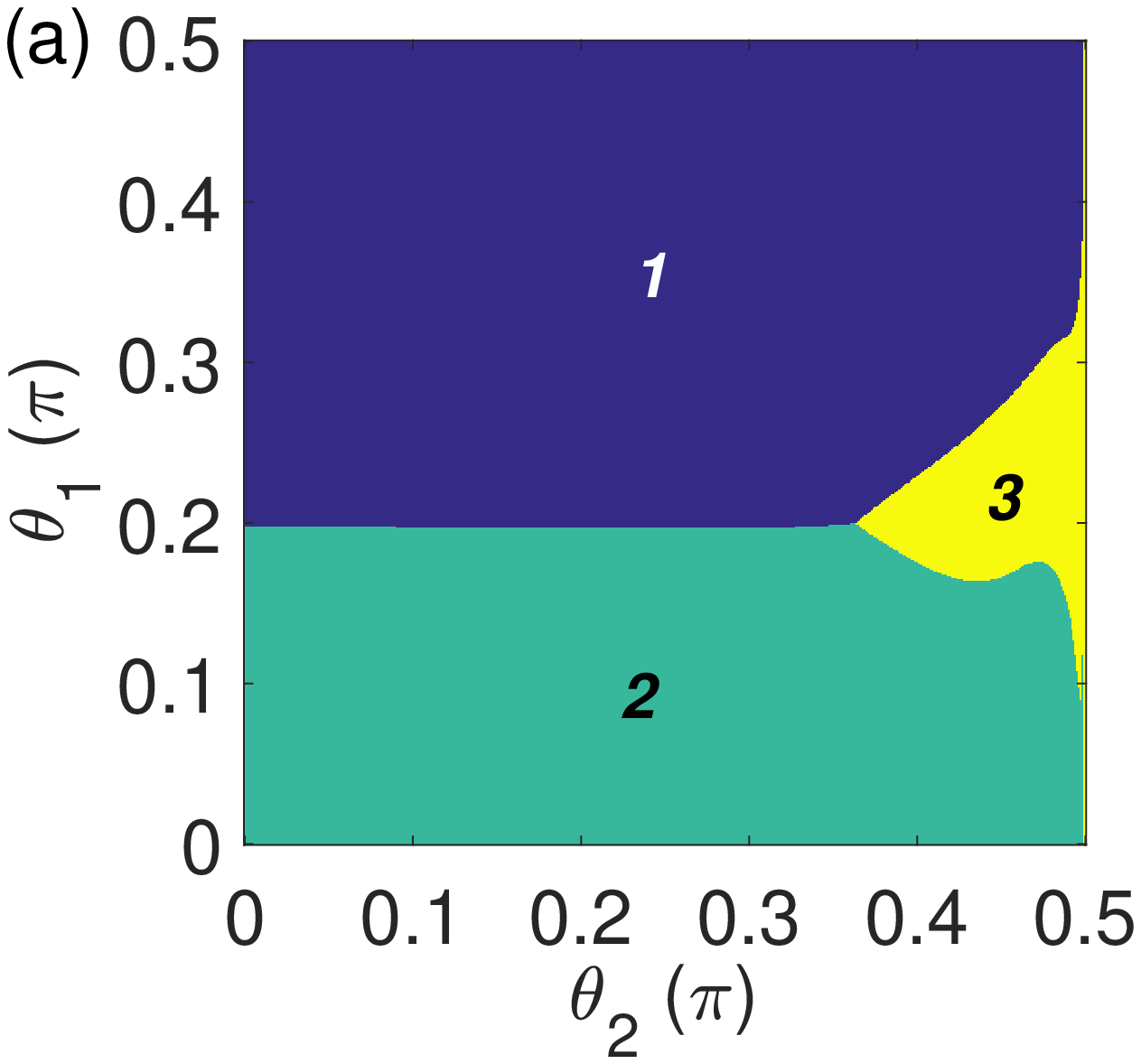}
\includegraphics*[width=4.2cm]{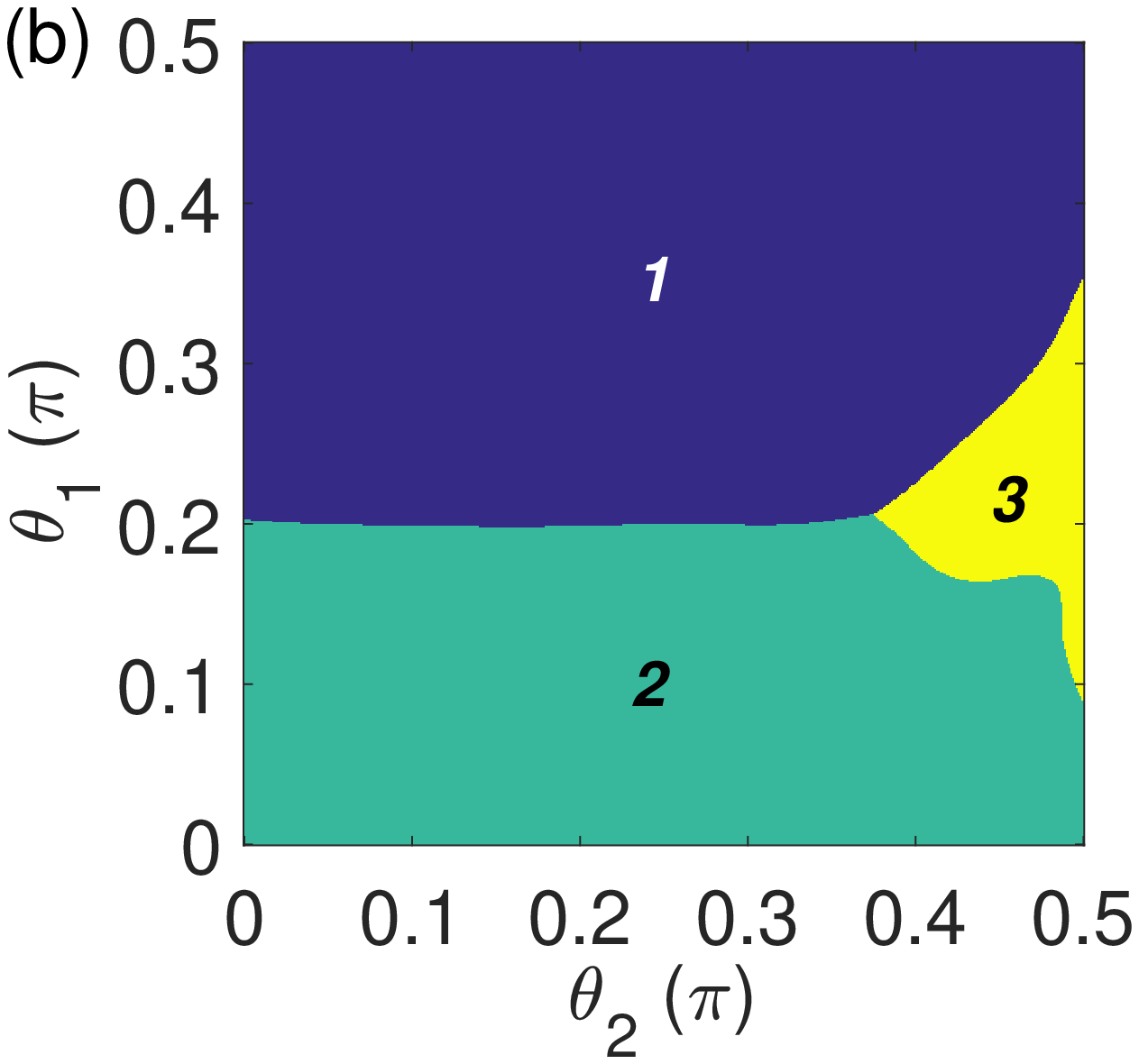}
\caption{The dependence of the control Hamiltonian choices on $\theta_1$ and $\theta_2$ with $\phi_1=\phi_2=0$ (blue for $H_1$, green for $H_2$ and yellow for low or the same fidelities). Results in (a) are calculated by simulating the dynamics with $H_1$ and $H_2$ and results in (b) are predicted by the neural network. Each subfigure contains $500\times500$ pixels corresponding to different $\theta_1$ and $\theta_2$.}
\label{FIG:Comparison}
\end{figure}

For illustrating the training process, we generated a training set with $10^4$ samples (59\% $H_1$ , 37\%$H_2$ and 4\% low or the same fidelities) plus a testing set with $N_{test}=5\times10^3$ samples through simulating the dynamics. In our simulations, the control parameters are  $\omega_2=2\omega_1$, $\omega_3=5\omega_1$, $g=0.5\omega_1$, $K=1$, and the control time is $T=20/\omega_1$. We set up a feedforward neural network with 4 input nodes, 3 output neurons and 2 hidden layers of 30 neurons, respectively. The feedforward neural network was trained by a back  propagation algorithm to minimize the MSE for the training set where a gradient descent method with momentum and adaptive learning rate was used. The training process is illustrated in Fig.\ref{FIG:Training} where the MSEs and the classification success rates $R$ for the training set and the testing set are plotted (with $10^4$ iterations). The MSEs for both set decreased dramatically in the first thousand iterations together with a rapid increase of the classification success rates (exceed 97\%). The training set MSE monotonically decrease through the training process due to the gradient algorithm, whereas the testing set MSE might slightly oscillate in the later iterations.

\begin{table}[!ht]
 \caption{Feedforward neural network performances for different numbers of training samples}
 \centering
 \begin{tabular}{ccccc}
  \hline
   $N_{train}$ & $N_{test}$ & $MSE$ (Testing) & Iteration &$R_A$ \\
  \hline
  $1\times10^1$        & $5\times10^3$ & 0.1060 & $1.5\times10^3$ & 78.3\% \\
  $1\times10^2$        & $5\times10^3$ & 0.0490 & $3.5\times10^4$ & 90.5\%  \\
  $1\times10^3$        & $5\times10^3$ & 0.0171 & $2.4\times10^3$ &  96.8\%  \\
  $1\times10^4$        & $5\times10^3$ & 0.0051 & $1.92\times10^4$ &98.7\%   \\
  $4\times10^4$ & $5\times10^3$ & 0.0034 & $1\times10^5$ &99.3\%  \\
  \hline
 \end{tabular}
\label{tab:classtests}
\end{table}

We then conducted 5 studies with different numbers of training samples $N_{train}$ where the testing sample numbers $N_{test}=5\times10^3$ are the same for comparison. The control parameters and the neural network structure are the same as those in Fig.\ref{FIG:Training}. In these studies, we adopted the iteration numbers corresponding to the minimal testing set MSEs (in at most $1\times10^5$ iterations) to determine the weights and biases of the neural networks. Finally, the trained neural networks were applied to predict control Hamiltonians for $5\times10^4$ new random initial states as applications. The corresponding classification success rates, denoted by $R_A$, and other training details are shown in Table.\ref{tab:classtests}. It is seen that the success rate $R_A$ was greater than $90\%$ even with 100 training samples in this example. In these studies, the MSEs for the testing set almost decreased to their minimums with thousands of iterations. We then  checked the dependence of the control Hamiltonian choice on $\theta_1$ and $\theta_2$ ($\phi_1=\phi_2=0$ is set) predicted by the neural network trained with $4\times10^4$ samples. The result is similar to the real result (calculated by simulating the dynamics) as shown in Fig.\ref{FIG:Comparison}. The processing time of the feedforward neural network depends on its layer number and node number in each layer. For Fig.\ref{FIG:Comparison}, the processing time of the feedforward neural network for Fig.\ref{FIG:Comparison}(b) was typically 1.5-5.5 orders lower than the simulation time for Fig.\ref{FIG:Comparison}(a) on our computer, depending on how the input vectors were sent to the neural network function, one-by-one or in batch.

\subsection{Designing Lyapunov functions}
In the second illustration, we will use GRNN to design initial-state-adaptive Lyapunov functions where one control Hamiltonian $H_1$ given by Eq.(\ref{eqn:Hami1}) is used. The control parameters are $\omega_2=2\omega_1$, $\omega_3=5\omega_1$, $g=0.5\omega_1$, $K=1$ and $T=2\pi/\omega_1$. The operator $P$ is designed by $P=p_1\ket{E_1}\bra{E_1}+p_2\ket{E_2}\bra{E_2}$ ($p_3=0$) according to Eq.(\ref{eqn:Pdesign}).

We generated totally $10^5$ samples for training where the input vector is $X^k=[\theta^k_{1}\ \theta^k_{2}\ \phi^k_{1}\ \phi^k_{2}]^T$ and the target output vector is $Y^k=[p^k_1\ p^k_2]^T$. Here $p_1^k$ and $p_2^k$ correspond to $p_1$ and $p_2$, respectively, which were found by minimizing the infidelity with the interior-point algorithm using the MATLAB optimization toolbox. For each random initial state, 8 optimizations with different (random) starting points were implemented to avoid local minimums. $p_{1,2}$ were optimized with the constraints $0\leq p_1\leq10, 0\leq p_2\leq20$. It is observed that the optimal values of $p_{1,2}$ corresponding to over $96\%$ of the initial states lie inside the area given by the constraints (rather than near the edges of the area), implying the constrains are plausible. The averaged logarithmic infidelity $\epsilon$ (defined by Eq.(\ref{eqn:AveLogInf})) for these training samples is $-4.03$ and the fraction of fidelities that are greater than $0.999$ is $R_{F>0.999}=0.773$.

\begin{figure}
\includegraphics*[width=7cm]{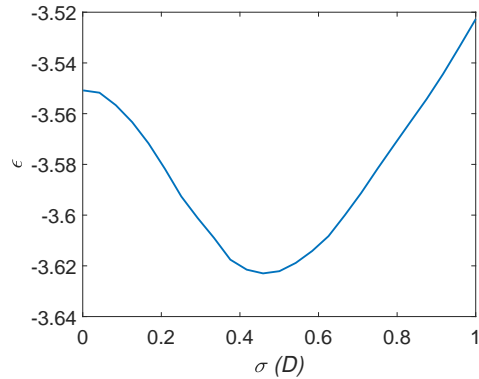}
\caption{The dependence of the averaged logarithmic infidelity $\epsilon$ for $2\times10^3$ testing samples on the smoothing parameter $\sigma$ (defined in Eq.(\ref{eqn:GRNN1})) of the GRNN.}
\label{FIG:GRNNtrain}
\end{figure}

\begin{figure}
\includegraphics*[width=9cm]{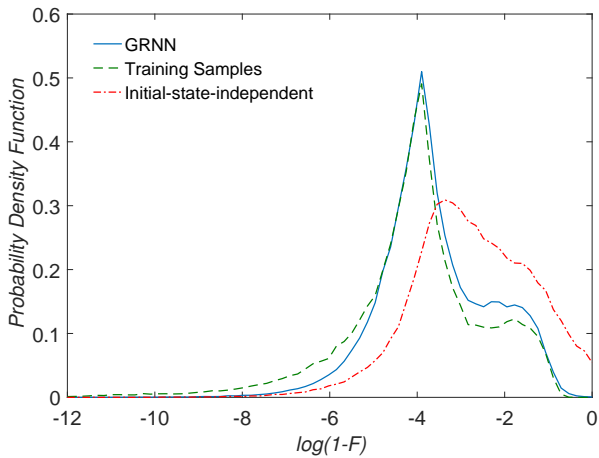}
\caption{Infidelity distributions for $10^5$ random initial states with different control schemes. The blue (solid) line and green (dashed) line represent the distribution from the GRNN control and that of the training samples (by numerical optimizations), respectively. The red (dot-dashed) line represents the result with an initial-state-independent (optimized) Lyapunov function. }
\label{FIG:dist}
\end{figure}

With a training set of $N_{train}=5\times10^4$ samples, we set up a GRNN with 4 input nodes, 2 output neurons, 3 summation layer neurons and $N_{train}$ pattern layer neurons. The smoothing parameter $\sigma$ was determined by checking the averaged logarithmic infidelity $\epsilon$ for $N_{test}=2\times10^3$ random initial states (out of the training set) with several trials and choosing the smoothing parameter with the minimal $\epsilon$. The dependence of $\epsilon$ on the smoothing parameter $\sigma$ is shown in Fig.\ref{FIG:GRNNtrain} with $0.001D\leq\sigma\leq D$ where $D=2/\sqrt[4]{N_{train}}$ is the average space between two neighboring normalized input vectors intruduced in Sec.\rm{IV} B. Finally, $\sigma= 0.46D$ was adopted for the GRNN. As an application, the trained GRNN was used to predict Lyapunov functions for $N_{app}=10^5$ new (random) initial states. The averaged logarithmic infidelity  and the fraction of $F>0.999$ is $\epsilon=-3.61$ and $R_{F>0.999}=0.715$, respectively.  To further demonstrate the performance of the GRNN, we plot the distribution of the infidelities of the application procedure in Fig.\ref{FIG:dist}. This distribution is compared with that of the  $10^5$  training samples and that from an initial-state-independent Lyapunov control with $P_{ind}=0.759\ket{E_1}\bra{E_1}+3.683\ket{E_2}\bra{E_2}$  for $N_{app}=10^5$ random initial states (see Fig.\ref{FIG:dist}). Here $P_{ind}$ was obtained by a numerical optimization that minimized the averaged logarithmic infidelity for $2\times10^3$ testing random initial states. The infidelity distribution from the GRNN-designed control is similar to that of the training samples, both with a peak at about $\log(1-F)=-4$, while the probability density functions with $\log(1-F)<-5$ and $-3<\log(1-F)<-1$ are slightly different. In contrast, the initial-state-independent control generated more low-fidelity states with an averaged logarithmic infidelity $\epsilon=-2.76$ and $R_{F>0.999}=0.443$, although $P_{ind}$ had been optimized. We remind that an arbitrary unoptimized initial-state-independent $P$ leads to a much worse result. For example, $p_1=p_2=10$ lead to $\epsilon=-1.26$  and $R_{F>0.999}=0.048$ in a simulation with $N_{app}=10^5$ random initial states.

\begin{table}[!ht]
 \caption{GRNN performances for different numbers of training samples}
 \centering
 \begin{tabular}{cccccc}
  \hline
     $N_{train}$ &  $\sigma$  &   $\epsilon$   &  $R_{F>0.999}$  & $R_{F>0.999}(train)$  \\
  \hline
   $5\times10^3$    & 0.50$D$  &  -3.41 & 0.657  &  0.773 \\
   $1\times10^4$   &  0.50$D$ &  -3.47 & 0.678  &  0.773 \\
   $5\times10^4$   &  0.50$D$ &  -3.60 & 0.716  &  0.773 \\
   $1\times10^5$   &  0.50$D$ &  -3.64 & 0.722  &  0.773 \\
  \hline
 \end{tabular}
\label{tab:GRNNtests}
\end{table}

In this example, we found empirically that the optimal smoothing parameters are generally near $0.5D$ either use $\epsilon$ or the MSE as the performance measure, regardless of the training sample number. The reason is that the distribution of the (normalized) training input vectors in the GRNN is known (with an average space $D$) and when $\sigma=0.5D$, the full width at half maximum of the Gaussian function in Eq.(\ref{eqn:GRNN1}) is roughly $D$. Such a Gaussian function probably leads to a good performance of a GRNN which could be understood by examining curve fitting problems. Thus one might find optimal $\sigma$ near $0.5D$ or use $0.5D$ for simplicity.

We further tested the performances of several GRNNs based on different numbers of training samples with $\sigma=0.5D$. The GRNNs were applied to the same Lyapunov control problem for $N_{app}=10^5$ random initial states. The details are shown in Table \ref{tab:GRNNtests}. As the number of training samples increases, the performance of the GRNN became better, together with a longer processing time of the GRNN due to its increased size. In our study, the longest prediction time (for one initial state) of the GRNN based on $10^5$ training samples was roughly the time of simulating the dynamics once. Meanwhile, finding an optimal $P$ in our  numerical optimizations (with 8 starting points) typically required simulating the dynamics about 6 hundred times. Thus, the initial-state-adaptive control with GRNN is able to improve the control performance with significantly less computing resources compared with that by numerical optimizations.

\section{Summary and discussion}
We have proposed two initial-state-adaptive Lyapunov control designs with machine learning where the feedforward neural network and the GRNN are used to select control schemes and predict Lyapunov functions for different initial states. The aim of the designs is to improve the control performance for different initial states without much increase of computing resources. Our methods could be applied to Lyapunov control problems when many initial states are involved or real-time processing is needed. The neural networks are trained by samples which are numerically generated before the final applications. We illustrated our designs with a three-level eigenstate control problem. Our results show that the neural networks are able to effectively learn the relationship between the initial states and the optimal control schemes or optimal Lyapunov functions and make predictions. The processing time of the neural networks are significantly less than that by numerical methods in our examples.

In our examples, the samples were divided into one training set and one testing set where we have assumed that the number of testing samples is large enough to reflect the generalization ability of the neural networks. When the number of samples is limited, other methods such as a k-fold cross validation \cite{Alpaydin2010} might  be used to fully take advantage of the samples. In our examples, the raw training data generated from simulations were used for the GRNN. In fact, our investigations showed that some simple processing of the raw training data may further improve the GRNN performance, for example, modestly removing the samples near the edge of the search area (given by $0<p_l\leq p_l^{max}$). The reason is that the relation between the optimal parameters and initial states may become less noisy after the data processing, although the training sample number is reduced. In general, the number of training samples and the size of neural networks significantly increase with the system dimension, the problem might be circumvented by investigating the initial-state-adaptive control with a subspace of all the possible initial states. Our initial-state-adaptive Lyapunov designs might also be used for other Lyapunov control problems to classify the initial states or predict continues parameters, such as the operator in Lyapunov function and the control fidelity.

\section{ACKNOWLEDGMENTS}
This work is supported by the National Natural Science Foundation of China under Grant No. 11705026, 11534002, 11775048, 61475033, the China Postdoctoral Science Foundation under Grant No. 2017M611293, and the Fundamental Research Funds for the Central Universities under Grant No. 2412017QD003.

\end{document}